\documentclass{tran-l}

\usepackage{amsmath, amssymb}
\usepackage{graphicx}
\usepackage{color}
\usepackage{textcomp}
\usepackage{graphics}

 \newtheorem{thm}{Theorem}[section]

 \newtheorem{prop}[thm]{Proposition}
 \newtheorem{defn}[thm]{Definition}
 \newtheorem{rem}[thm]{Remark}
 \newtheorem{exa}[thm]{Example}
 \newcommand{\PP}{\mathbb{P}}

\newcommand{\bm}{\bibitem}
\newcommand{\no}{\noindent}
\newcommand{\be}{\begin{equation}}
\newcommand{\ee}{\end{equation}}

\newcommand{\bea}{\begin{eqnarray}}
\newcommand{\bes}{\begin{subequations}}
\newcommand{\ees}{\end{subequations}}
\newcommand{\bgt}{\begin{gather}}
\newcommand{\egt}{\begin{gather}}
\newcommand{\eea}{\end{eqnarray}}
\newcommand{\beaa}{\begin{eqnarray*}}
\newcommand{\eeaa}{\end{eqnarray*}}

\newcommand{\EE}{{\mathbb E}}
\newcommand{\RR}{{\mathbb R}}
\newcommand{\cal}{\mathcal}

\begin{document}

\title[(Martingale) Optimal Transport with neural networks]{(Martingale) Optimal Transport and anomaly detection with neural networks: a primal-dual algorithm}
\author{Pierre Henry-Labord\`ere}
\address{Soci\'et\'e G\'en\'erale, Global markets Quantitative Research}
\address{CMAP, Ecole Polytechnique}
\email{pierre.henry-labordere@sgcib.com} \subjclass{} 
\date{}

\maketitle

\begin{abstract}\no In this  paper, we introduce a  primal-dual algorithm for solving  (martingale) optimal transportation problems, with cost functions satisfying the twist condition, close to the one that has been used recently for training generative adversarial networks. As some additional applications, we consider anomaly detection and automatic generation of financial data.
\end{abstract}

\section{Introduction}

\no We introduce a  primal-dual algorithm for solving  (martingale) optimal transportation problem (in short MOT), potentially large-scale, using neural networks. The martingale optimal transport, first introduced  in \cite{bei} and in a continuous-time setting in \cite{gal},  can be defined in a discrete-time setting as the following infinite-dimensional linear program:
\bea \mathrm{MK}_c(\mu^1,\mu^2):=\sup_{\PP \in {\cal M}(\mu^1,\mu^2)} \EE^{\PP}[c(S_1,S_2)] \label{MOT} \eea where
$   {\cal M}(\mu^1,\mu^2):=\{ \PP \in {\cal P}(\RR^d,\RR^d) \; :\; S_1 \overset{\PP}{\sim} \mu^1, \quad S_2 \overset{\PP}{\sim}  \mu^2, \quad \EE^\PP[S_2|S_1]=S_1 \} $ is a weak compact convex set and  ${\cal P}(\RR^d \times \RR^d)$ is the set of probability measures on $\RR^d \times \RR^d$ (or $\RR_+^d \times \RR_+^d$ if the random variables $S_1$ and $S_2$ are interpreted as financial asset prices).  A similar definition applies by replacing the supremum over ${\cal M}(\mu^1,\mu^2)$ by an infimum.  $\mathrm{MK}_c(\mu^1,\mu^2)$ is a number which depends on a cost function $c: \RR^d \times \RR^d \mapsto \RR$ and two marginal distributions $\mu^1$ and $\mu^2$ defined on $\RR^d$. In comparison with the classical OT,  we have an additional martingale constraint $\EE^\PP[S_2|S_1]=S_1$ and the linear problem is well-posed if and only if $\mu^1 \leq \mu^2$ in the convex order. In mathematical finance, $\mathrm{MK}_c(\mu^1,\mu^2)$ can then be interpreted as the model-independent arbitrage-free optimal upper bound for a payoff $c(S_1,S_2)$ depending on an asset $S_\cdot \in \RR^d$ evaluated at two maturities $t_1<t_2$, i.e., $S_1:=S_{t_1}$, $S_2:=S_{t_2}$, which is  consistent with  the prices (at $t=0$) of $t_1$ and $t_2$ ($d$-dimensional) European basket options (see \cite{phl1} for an extensive introduction to MOT and its relevance in arbitrage-free pricing). Our algorithm, described in Section \ref{A primal-dual algorithm}, can also be  applied to more general linear programs of the form:
 \beaa \mathrm{P}_c:=\sup_{\PP \in {\cal M}} \EE^{\PP}[c(S_1,S_2,\cdots,S_n)] \label{P} \eeaa where ${\cal M}$ is a weak-compact convex subset of ${\cal P}((\RR^d)^n)$, see for example the multi-marginals (M)OT. However, our algorithm will be applicable only to cost functions satisfying a (martingale) twist condition. Although the extension of our algorithm to this more general setting is straightforward, we prefer for the sake of simplicity to focus on (martingale) OT as defined by (\ref{MOT}). Most of the numerical schemes of (M)OT, that we will describe, rely strongly on the dual Monge-Kantorovich formulation  in which $\mathrm{MK}_c(\mu^1,\mu^2)$ can be written as (see \cite{bei} for a proof in the context of MOT):
\bea \mathrm{MK}_c(\mu^1,\mu^2):=\inf_{u_1 \in \mathrm{L}^1(\mu^1),  u_2 \in \mathrm{L}^1(\mu^2), h \in C_b(\RR^d,\RR^d)}
\EE^{\mu^1}[u_1]+\EE^{\mu^2}[u_2] \label{MOTdual} \eea such that  for all $(s_1,s_2) \in \RR^d \times \RR^d$ \bea u_1(s_1)+u_2(s_2)+h(s_1).(s_2-s_1) \geq c(s_1,s_2) \label{linearconstraint} \eea \no By definition,  $h(s_1).(s_2-s_1):=\sum_{i=1}^d h_i(s_1)(s^i_2-s^i_1)$.

\section{Numerical algorithms: A short overview} \label{Numerical algorithms: A short overview}
\no In this section, we review three numerical algorithms for solving (martingale) optimal transport and highlight their main drawbacks\footnote{We acknowledge G. Peyr\'e for useful discussions.}. These algorithms will be compared to our primal-dual method in Section \ref{Numerical examples}.

\subsection{Simplex and cutting-plane} The problem (\ref{MOTdual}) (resp. \ref{MOT}) defines a linear program that can be solved using a simplex algorithm. In the context of MOT, this has been explored in \cite{phl2}. By discretizing the measures $\mu^1$ and  $\mu^2$ on a large grid $G_\infty$ in $\RR^d \times \RR^d$, we obtain a finite-dimensional linear program. Due to the large number $N:=\mathrm{card}(G_\infty)$ of linear constraints (\ref{linearconstraint}), one can use a cutting-plane algorithm, see \cite{phl2} for extensive details.  This  consists in solving  the LP program using first a small dimensional grid $G_0 \subset G_\infty$ ($\mathrm{card}(G_0) \ll \mathrm{card}(G_\infty)$). The optimal bound  $\mathrm{MK}^{(0)}_c(\mu^1,\mu^2)$ is attained by the  dual variables $(u_1^{(0)},u_2^{(0)},h^{(0)})$.  Then we check on the full grid $G_\infty$ if our optimal dual solution violates the linear constraints (\ref{linearconstraint}). The points of $G_\infty$ where the linear constraints  are not satisfied, are then added to the grid $G_0$, defining a new refined grid $G_1$.  By construction, we obtain $\mathrm{MK}_c(\mu^1,\mu^2)  \geq \mathrm{MK}^{(1)}_c(\mu^1,\mu^2) \geq \mathrm{MK}^{(0)}_c(\mu^1,\mu^2)$ as $G_0 \subset G_1 \subset G_\infty$. The procedure is then iterated until the optimal dual solution $(u_1^{(n)},u_2^{(n)},h^{(n)})$ at step $(n)$ satisfies all the constraints on $G_\infty$ for which  we can conclude that we have converged towards the true solution. Despite its simplicity, this algorithm could not be extended in large dimension as the number of constraints explodes with the dimension. For example, the complexity of the Hungarian/auction algorithms is $O(N^3)$.

\subsection{Entropic relaxation} \label{Entropy}
Another approach is to introduce an entropy penalization (or more generally a $f$-divergence):
\beaa \mathrm{MK}^\epsilon_c(\mu^1,\mu^2):=\sup_{\PP \in {\cal M}(\mu^1,\mu^2)} \EE^{\PP}[c(S_1,S_2)]- \epsilon H(\PP|\PP^0) \eeaa where
$H(\PP|\PP^0):= \EE^\PP[ \left( \ln {d \PP \over d\PP^0}-1 \right)]$ is the relative entropy with respect to a prior probability measure $\PP^0 \in {\cal P}(\RR^d \times \RR^d)$ and $\epsilon$ is a positive parameter taken to be small. In particular, $\lim_{\epsilon \rightarrow 0}\mathrm{MK}^\epsilon_c(\mu^1,\mu^2)=\mathrm{MK}_c(\mu^1,\mu^2)$.  The problem $\mathrm{MK}^\epsilon_c(\mu^1,\mu^2)$ can be dualized using the Fenchel-Rockafellar's theorem into a strictly convex optimization problem  \cite{phl2}:
\bea \mathrm{MK}^\epsilon_c(\mu^1,\mu^2):=\inf_{u_1 \in \mathrm{L}^1(\mu^1),  u_2 \in \mathrm{L}^1(\mu^2),  h \in C_b(\RR^d,\RR^d) }&&  \EE^{\mu^1}[u_1]+\EE^{\mu^2}[u_2] \nonumber \\&+&\epsilon  \EE^{\PP^0}[ e^{{1\over \epsilon}\left( c(s_1,s_2)-u_1(s_1)-u_2(s_2)-h(s_1).(s_2-s_1) \right)}] \label{ent} \eea

\subsubsection{Sinkhorn's algorithm}
By computing the gradients with respect to $u_1$, $u_2$ and $h$, we obtain the first-order optimality conditions:
\bea e^{-{u_1(s_1) \over \epsilon}} \int p_0(s_1,s_2)ds_2 e^{{ 1 \over \epsilon}\left( c(s_1,s_2)-u_2(s_2)-h(s_1).(s_2-s_1) \right)} &=& \mu^1(s_1) \label{u1} \\
e^{-{u_2(s_2) \over \epsilon}} \int p_0(s_1,s_2)ds_1 e^{{ 1 \over \epsilon}\left(c(s_1,s_2)-u_1(s_1)-h(s_1).(s_2-s_1)\right) } &=& \mu^2(s_2) \\
\int p_0(s_1,s_2)ds_2 (s_2-s_1)e^{{ 1 \over \epsilon}\left(c(s_1,s_2)-u_2(s_2)-h(s_1).(s_2-s_1)\right) } &=&0
 \eea   For the sake of simplicity, we have assumed here  that $\PP^0$, $\mu^1$ and $\mu^2$ are absolutely-continuous with respect to the Lebesgue measure. The Sinkhorn algorithm can be then described by the following steps:
 \begin{enumerate}
 \item Set $n:=1$ and set $u^{(0)}_1:=0$, $u^{(0)}_2:=0$, $h^{(0)}:=0$ for convenience. We approximate the measures $\mu^1$ and $\mu^2$ by Dirac masses supported on $N$ points $(s_1^i)_{1 \leq i \leq N}$ and $(s_2^i)_{1 \leq i \leq N}$.
 \item Compute $u_1^{(n)}(s_1)$ for all $(s_1^i)_{1 \leq i \leq N}$  using
 \beaa e^{-{u^{(n)}_1(s_1) \over \epsilon}} \int p_0(s_1,s_2)ds_2 e^{{ 1 \over \epsilon}\left( c(s_1,s_2)-u^{(n-1)}_2(s_2)-h^{(n-1)}(s_1).(s_2-s_1) \right)} &=& \mu^1(s_1) \eeaa
  \item Compute $h^{(n)}(s_1)$  for all $(s_1^i)_{1 \leq i \leq N}$ by finding the (unique) zero $\theta \in \RR^d$ of
 \beaa h(s_1):=\theta \; \mathrm{s.t.} \; \int p_0(s_1,s_2)ds_2 (s_2-s_1)e^{{ 1 \over \epsilon}\left(c(s_1,s_2)-u^{(n-1)}_2(s_2)-\theta.(s_2-s_1)\right) } &=&0\eeaa
 \item Compute $u_2^{(n)}(s_2)$ for all $(s_2^i)_{1 \leq i \leq N}$ using
 \beaa e^{-{u^{(n)}_2(s_2) \over \epsilon}} \int p_0(s_1,s_2)ds_1 e^{{ 1 \over \epsilon}\left(c(s_1,s_2)-u^{(n)}_1(s_1)-h^{(n)}(s_1).(s_2-s_1)\right) } &=& \mu^2(s_2) \eeaa
\item Set $n:=n+1$ and iterate steps (2-3-4) up to convergence.
  \end{enumerate} The use of the Sinkhorn algorithm for solving OT problem was introduced in \cite{cut} and in \cite{guo}, \cite{demarch0} in the context of MOT (see also \cite{demarch1} for an application to the construction of arbitrage-free implied volatility surfaces). Again this algorithm does not scale well with the dimension as at each Sinkhorn's iteration,
  $u_1^{(n)}(s_1), h^{(n)}(s_1), u_2^{(n)}(s_1)$ must be computed on a grid whose the cardinality explodes with the dimension $d$. The overall complexity is $O(N^2 \ln N)$.

\subsection{and neural networks...}

\no In \cite{seg}, the optimization (\ref{ent}) is solved by approximating the potentials $u_1,u_2$ (and $h$) by some neural networks and then the training is achieved  using a stochastic gradient descent algorithm. Similarly, by using Equation (\ref{u1}),   the problem (\ref{ent}) can be converted into an equivalent form which involves only the potentials $u_2$ and $h$:
\beaa \mathrm{MK}^\epsilon_c(\mu^1,\mu^2)&:=&\inf_{ h \in C_b(\RR), u_2 \in \mathrm{L}^1(\mu^2)} \EE^{\mu^2}[u_2] \\&&\epsilon \int \mu^1(ds_1)
\ln \int p^0(s_1,s_2)ds_2 e^{{1 \over \epsilon}\left( c(s_1,s_2)-u_2(s_2)-h(s_1).(s_2-s_1) \right)} \\
&-&\epsilon \int \mu_1(ds_1)\left( \ln \mu_1(s_1)-1 \right) \eeaa and solve similarly. In \cite{gen}, instead of using neural networks, the authors make use of an expansion of the dual variables in a reproducing
kernel Hilbert space. Despite this algorithm scales properly with the dimension in practise,  we will illustrate  in our numerical experiments that our computations are unstable when $\epsilon$ becomes small. This has been also reported in \cite{gen}.

\subsection{Penalization} In \cite{eck}, the optimization $\mathrm{MK}_c(\mu^1,\mu^2)$ is approximated by
\beaa \mathrm{MK}^\gamma_c(\mu^1,\mu^2):&=&\inf_{u_1 \in \mathrm{L}^1(\mu^1),  u_2 \in \mathrm{L}^1(\mu^2), h \in C_b(\RR)} \EE^{\mu^1}[u_1]+\EE^{\mu^2}[u_2]\\
&+&\gamma \EE^{\PP^0}[ (c(s_1,s_2)-u_1(s_1)-u_2(s_2)-h(s_1)(s_2-s_1))_+^2 ] \eeaa where $\gamma$ is a large parameter. This ensures that by taking $\gamma$ large, the optimal dual solution $(u_1^*,u_2^*,h^*)$ will satisfy the linear constraints (\ref{linearconstraint}) and therefore $\lim_{\gamma \rightarrow \infty} \mathrm{MK}^\gamma_c(\mu^1,\mu^2)= \mathrm{MK}_c(\mu^1,\mu^2)$. As above, the potentials $u_1$, $u_2$ and $h$ are approximated by  some neural networks. This is a classical technique for solving linear programs by penalization and in practise the parameter $\gamma_t$ is chosen to increase to a large value as the learning parameter $\eta_t$, used in the stochastic gradient descent, decreases.  In our numerical experiments,
we will illustrate that this algorithm is unstable, when the parameter $\gamma$ is chosen large in order to converge to the true solution. Finally, let us remark that the
penalization method can be obtained by replacing the entropy penalization $H(\PP|\PP^0)$ by the $\mathrm{L}^2$-divergence
$ f(\PP|\PP^0):=\EE ^{\PP^0}[ \left({d\PP \over d\PP^0}\right)^2]  $.

\label{Penalization}

\section{A primal-dual algorithm} \label{A primal-dual algorithm}
\subsection{A saddle-point formulation}
\no For the sake of clarity, we  explain our algorithm in the case of the classical OT problem which consists in solving
\beaa \mathrm{MK}_c(\mu^1,\mu^2):=\sup_{\PP \in {\cal M}(\mu^1,\mu^2)} \EE^{\PP}[c(S_1,S_2)] \eeaa where
$  {\cal M}(\mu^1,\mu^2):=\{ \PP \in {\cal P}(\RR^d,\RR^d) \; :\; S_1 \overset{\PP}{\sim} \mu^1, \quad S_2 \overset{\PP}{\sim}  \mu^2 \}$.  By introducing the Lagrange multipliers $u_1$ and $u_2$ associated to the two marginal constraints,  this problem can be written as a minimax (relaxed) optimization problem:
\bea \mathrm{MK}_c(\mu^1,\mu^2):&=&
\inf_{u_1 \in \mathrm{L}^1(\mu^1),  u_2 \in \mathrm{L}^1(\mu^2)} \sup_{\PP \in {\cal M}_+}
\EE^{\mu^1}[u_1]+\EE^{\mu^2}[u_2] \nonumber \\&+&\EE^{\PP}[c(S_1,S_2)-u_1(S_1)-u_2(S_2)] \label{sp} \eea
\no where ${\cal M}_+$ denotes the space of positive measures on $\RR^d \times \RR^d$.

\subsection{Using Brenier's theorem}

\begin{defn}[Twist condition] A function $c \in C(\RR^d \times \RR^d)$ differentiable with
respect to $s_1$ is said to be twisted if $ \forall s_0  \in \RR^d$, the map $s_2 \in \RR^d \mapsto
 \nabla _{s_1}c(s_0, s_2)$ is one-to-one.
\end{defn}
\no We recall the Brenier theorem (see e.g. \cite{vil}):
\begin {thm}[Brenier's theorem]
\no By assuming that $\mu^1$ is absolutely continuous with respect to the Lebesgue measure and the cost function $c$ satisfies the twist condition, the optimal probability measure $\PP^*$, solution of the above saddle-point  problem (\ref{sp}), is supported on a  unique map $T: \RR^d \mapsto \RR^d$:
\beaa \PP^*(ds_1,ds_2)=\mu^1(ds_1)\delta(s_2-T(s_1))ds_2 \eeaa
\end{thm}
\no Note that  the constraints $S_1 \overset{\PP^*}{\sim} \mu^1$ and  $S_2 \overset{\PP^*}{\sim} \mu^2$  imply the requirement $T_\#\mu^1=\mu^2$ where $T_\#\mu^1$ denotes the push-forward of the measure $\mu_1$ by the map $T$. $T$ can be characterized as the unique solution of a Monge-Amp\`ere-like equation. More precisely, in the case of the quadratic cost function, $T$ is the gradient of a convex  function solution of the Monge-Amp\`ere PDE (see e.g. \cite{vil}).

\begin{rem}[Fr\'echet-Hoeffding $d=1$] Under the (twist) condition $\partial_{s_1 s_2}c \geq 0$ in $d=1$, the optimal transport can be solved analytically and it is given by the Fr\'echet-Hoeffding solution:
\bea \mathrm{MK}_c(\mu^1,\mu^2)=\int_0^1 (F_1^{-1}(u)-F_2^{-1}(u))^2 du \label{HF} \eea The map is then $T(s)=F_2^{-1} \circ F_1(s)$ with $F_i$ the cumulative distribution of $\mu^i$.
\end{rem}
\no Under the twist condition, the above minimax optimization (\ref{sp}) can therefore be simplified as
\bea \mathrm{MK}_c(\mu^1,\mu^2):=\inf_{u \in \mathrm{L}^1(\mu^2)} \sup_{T: \RR^d \mapsto \RR^d}  \EE^{\mu^1}[ c(S_1,T(S_1))-u(T(S_1)) ]+\EE^{\mu^2}[u(S_2)] \label{brenier} \eea Note that as $S_1 \overset{\PP^*}{\sim} \mu^1$, the potential $u_1$ has disappeared and the minimax optimization involves now only the potential $u:=u_2$ and the Brenier map $T$.

\subsection{and neural networks...}

We then approximate the two unknowns $u:\RR^d \mapsto \RR$ and $T:\RR^d \mapsto \RR^d$ with two neural networks depending respectively on some weights $\theta \in \RR^u$ and $\omega \in \RR^t$.  $\mathrm{MK}_c(\mu^1,\mu^2)$ can then be approximated by
\bea \mathrm{MK}^{t,u}_c(\mu^1,\mu^2):=\min_{\theta \in \RR^u} \max_{\omega \in \RR^t}  \EE^{\mu^1}[ c(S_1,T_{\omega}(S_1))-u_{\theta}(T_\omega(S_1)) ]+\EE^{\mu^2}[u_{\theta}(S_2)] \label{nn} \eea In particular, from the universal approximation property of neural networks, we have
$\lim_{t,u \rightarrow \infty}\mathrm{MK}^{t,u}_c=\mathrm{MK}_c$.

\subsection{Link with Wasserstein generative adversarial networks}

\no The $p$-Wasserstein distance ${\cal W}_p(\mu^1,\mu^2) $ corresponds to an OT problem with a $\mathrm{L}^p$-cost in $\RR^d$, $c(s_1,s_2):=|s_2-s_1|^p$: \beaa \left({\cal W}_p(\mu^1,\mu^2)  \right)^p:=\inf_{\PP \in {\cal M}(\mu^1,\mu^2)} \EE^{\PP}[|S_2-S_1|^p] \eeaa ${\cal W}_p$ defines then a distance which metrizes the space  ${\cal P}(\RR^d)$ (see e.g. \cite{vil}). If we consider a probability measure $\mu^\mathrm{real}$ in $\RR^d$ corresponding to some real data, one would like to reconstruct this density using a mapping $\hat{T}: \RR^l \mapsto \RR^d$ with $l \ll d$ and such that the push-forward of $\hat{T}$ by a prior density $\mu^0$ supported on $\RR^l$  (e.g. an uniform or Gaussian density for the sake of simplicity) is as close as possible to $\mu^\mathrm{real}$ with respect to the Wasserstein distance. The mapping $\hat{T}$ is then chosen to be the solution  of
\beaa \mathrm{P}&:=&\inf_{\hat{T}: \RR^l \mapsto \RR^d} {\cal W}_p(\mu^1,\mu^2)  \\
&=& \inf_{\hat{T}: \RR^l \mapsto \RR^d}  \inf_{\PP \in {\cal M}(\mu^\mathrm{real},\hat{T}_\# \mu^0)} \EE^{\PP}[|S_2-S_1|^p] \eeaa Note that $H(\hat{T}_\# \mu^0|\mu^\mathrm{real})=+\infty$ and this is why it is not possible to use the relative entropy as  in the case of maximum likelihood estimation.
\no Using the saddle-point formulation of the Wassertein distance (the $\mathrm{L}^p$-cost satisfies the twist condition) explained in the previous section, this is equivalent to the following minimax optimization:
\beaa \mathrm{P}=\sup_{u\in \mathrm{L}^1(\mu^\mathrm{real})} \inf_{\hat{T}: \RR^l \mapsto \RR^d,T: \RR^d \mapsto \RR^d}   \EE^{\mu^\mathrm{real}}[ c(S_1,T(S_1))-u(T(S_1))]  +\EE^{\mu^0}[u(\hat{T}(S_0))] \eeaa This problem is similar to (\ref{brenier}) and therefore as described in Section \ref{Arrow-Hurwicz algorithm: recipe}, our algorithm is close in spirit to the one used for training Wasserstein generative adversarial networks \cite{arj} (see also \cite{goo}). %Using neural network approximations for $T,\hat{T}$ and $u$, we get
%\bea \mathrm{P}^{t,\hat{t},u}:=\max_{\theta \in \RR^u} \min_{\omega \in \RR^t, \hat{\omega} \in \RR^{\hat{t}} }  
%\EE^{\mu^\mathrm{real}}[ c(S_1,T_{\omega}(S_1))-u_{\theta}(T_\omega(S_1)) ]+%\EE^{\mu^0}[u_{\theta}(\hat{T}_{\hat{\omega}}(S_0))] \eea

\no  Specializing to $p=1$, we get
\beaa \mathrm{P}=\sup_{u\in \mathrm{L}^1(\mu^\mathrm{real})} \inf_{\hat{T}: \RR^l \mapsto \RR^d,T: \RR^d \mapsto \RR^d}   \EE^{\mu^\mathrm{real}}[ |S_1-T(S_1)|-u(T(S_1))]  +\EE^{\mu^0}[u(\hat{T}(S_0))] \eeaa This should be compared with the dual formulation of the $1$-Wassertein distance used in \cite{arj}
\beaa \mathrm{P}=\sup_{u\in \mathrm{Lip}_1} \inf_{\hat{T}: \RR^l \mapsto \RR^d}   -\EE^{\mu^\mathrm{real}}[ u(S_1)]  +\EE^{\mu^0}[u(\hat{T}(S_0))] \eeaa where the supremum is over all the $1$-Lipschitz functions.  The Lipschitz constraint is enforced in brute force by weight clipping.

\no Starting from the primal formula of OT and using the Brenier theorem, $\mathrm{P}$ can also be written as
 \beaa \mathrm{P}&=&  \inf_{\hat{T}: \RR^l \mapsto \RR^d,T: \RR^d \mapsto \RR^d s.t. T_\#\mu^\mathrm{real}=\hat{T}\#\mu^0}   \EE^{\mu^\mathrm{real}}[ c(S_1,T(S_1))]    \eeaa This was done in \cite{bou} although the  Brenier result is not mentioned.  The constraint $T_\#\mu^\mathrm{real}=\hat{T}_\#\mu^0$ is then implemented by adding a penalty term $\gamma D(\cdot|\cdot)$ with $\gamma$ large:
 \beaa \mathrm{P}^\gamma&:=&  \inf_{\hat{T}: \RR^l \mapsto \RR^d,T: \RR^d \mapsto \RR^d }   \EE^{\mu^\mathrm{real}}[ c(S_1,T(S_1))]   +\gamma
 D(T_\#\mu^\mathrm{real}|\hat{T}_\#\mu^0) \eeaa One  obtains the Wasserstein-VAE formulation.

\subsection{Anomaly detector and data generator} Let us consider some real data generated by a density $\mu^\mathrm{real}$ and let us choose a prior density $\mu^0$ supported on a low-dimensional manifold. As outlined above, we find the density $\hat{T}_\# \mu^0$ such that the $p$-Wasserstein distance ${\cal W}_p(\mu^\mathrm{real},\hat{T}_\# \mu^0)$ is minimized.  Then, a data $x_\mathrm{anomaly}$ will be considered as an anomaly  if $\hat{T}_\# \mu^0(x_\mathrm{anomaly})$ is below a certain threshold $\lambda$:
\beaa \hat{T}_\# \mu^0(x_\mathrm{anomaly}) \leq \lambda \eeaa
Similarly, a new data $x_\mathrm{new}$ can be generated by  drawing a random variable $Z$ distributed according to $\mu^0$ and set $x_\mathrm{new}=\hat{T}(Z)$.  \label{Anomaly detection and data generator}

\subsection{Arrow-Hurwicz algorithm: recipe} \label{Arrow-Hurwicz algorithm: recipe}

\no We simulate $\mu^1$ and $\mu^2$ by Monte-Carlo with $N_\mathrm{MC}$ paths $(S^i_1,S^i_2)_{1 \leq i \leq N_\mathrm{MC}}$ and for large $N_\mathrm{MC}$, our optimization (\ref{nn}) consists in solving:
\beaa \min_{\theta \in \RR^u} \max_{\omega \in \RR^t}   {1 \over N_\mathrm{MC}}\sum_{i=1}^{N_\mathrm{MC}} J_i(\theta,\omega)\eeaa
\no where
\beaa  J_i(\theta,\omega):=  c(S^i_1,T_\omega(S^i_1))-u_\theta(T_\omega(S^i_1)) + u_\theta(S^i_2)  \eeaa The average functional  can be optimized by using a stochastic Arrow-Hurwicz algorithm which consists in doing sequentially the two iterations at each step $n$: Draw a uniform r.v. $I \in [[1,N_\mathrm{MC}]]$ and compute
\bea \theta_{n+1}&=&\theta_n -\eta \nabla_\theta J_I(\theta_{n},\omega_{n}) \label{AH1}\\
\omega_{n+1}&=&\omega_n +\eta\nabla_\omega J_I(\theta_{n+1},\omega_{n}) \label{AH2}
\eea where $\eta$ is a  learning parameter.  In practise, the gradients are computed by back-propagation where
\beaa  \nabla_\theta J_I(\theta,\omega)&=&  - \nabla_\theta u_\theta(T_\omega(S^I_1)) +  \nabla_\theta u_\theta(S^I_2) \\
 \nabla_\omega J_I(\theta,\omega)&=& \left( \nabla_{s_2} c(S^I_1,T_\omega(S^I_1))
  -\nabla_{s_2} u_\theta(T_\omega(S^I_1))\right).\nabla_\omega T_\omega(S^I_1 )  \eeaa

\no We could  used also a predictor-corrector scheme (that gives similar results in our numerical experiments):
\beaa \theta_{n+1/2}&=&\theta_n -\eta  \nabla_\theta J_I(\theta_{n},\omega_{n}) \\
\theta_{n+1}&=&\theta_n -\eta  \nabla_\theta J_I(\theta_{n+1/2},\omega_{n}) \\
\omega_{n+1/2}&=&\omega_n +\eta  \nabla_\omega J_I(\theta_{n+1},\omega_{n})  \\
\omega_{n+1}&=&\omega_n +\eta \nabla_\omega J_I(\theta_{n+1},\omega_{n+1/2})
\eeaa

\subsection{Convergence}
\no By using one layer for the approximation of the two unknowns $T_\omega$ and $u_\theta$ with a linear activation function (a drift can also be  included without loss of generality):
\beaa T(x)&:=&\omega.  x,\quad  u(x):=\theta^\dag. x, \quad \omega \in \mathrm{M}_{p,p}, \quad  \theta \in \RR^p
\eeaa the problem (\ref{nn}) can be written as
\bea  \min_{\theta \in \RR^p} \max_{\omega \in \mathrm{M}_{p,p}}  \EE^{\mu^1}[c(X,\omega X)]-\theta^\dag \omega  \EE^{\mu^1}[X] +\theta^\dag \EE^{\mu^2}[X] \label{cha} \eea and it is of the form
\beaa \min_x \max_y  y^\dag K x +G(x)+F(y) \eeaa where $K$ is a linear operator. As shown by \cite{cha}, the stochastic Arrow-Hurwicz algorithm converges if $F$ is concave, $G$ is convex and $||K||\eta^2<1$. Our program (\ref{cha}) is clearly convex in
$\theta$ as being linear and is concave in $\omega$ if and only if $D^2 _{s_2} c \leq 0$. This implies that our algorithm  converges (in the case of one layer), if we impose that  $D^2 _{s_2} c \leq 0$. Additionally, we should have that $c$ satisfies the twist condition  as we have used the Brenier theorem.

\no Let us remark that if we consider the new cost function $\bar{c}(s_1,s_2)=c(s_1,s_2)-U(s_2)$, then we have for all $U \in \mathrm{L}^1(\mu^2)$:
\beaa \mathrm{MK}_{\bar{c}}(\mu^1,\mu^2)+\EE^{\mu^2}[U(S_2)]=\mathrm{MK}_c(\mu^1,\mu^2) \eeaa
Using this property, we can apply our algorithm to the cost function $\bar{c}$ where $U$ is chosen such that\footnote{We are grateful to our master-degree students Y. Chen and F. Jiang at Ecole Polytechnique for pointing to us this remark.}
\beaa D^2 _{s_2}\bar{c}=D^2 _{s_2} c-D^2 _{s_2} U(s_2) \leq 0,\quad \forall \;  (s_1,s_2) \in \RR^d \times \RR^d \eeaa

\begin{exa} For $c(x,y)=-(x-y)^2$, we can take $U(y)=0$. For $c(x,y)=(x+y)^2$, we can take $U(y)=2 y^2$.
\end{exa}

\no Using the result in \cite{cha}, we conclude:
\begin{prop}[Convergence] Let us assume that $c$ satisfies the twist condition and $D^2 _{s_2} c-D^2 _{s_2} U(s_2) \leq 0$  for some twice differentiable function $U$ in $\mathrm{L}^1(\mu^2)$,  then the Arrow-Hurwicz algorithm (with one layer) (\ref{AH1}-\ref{AH2}) converges for $\eta$ small enough. \label{Convergence}
\end{prop} \no Note that a similar conclusion appears if we expand $T_\omega$ and $u_\theta$ in terms of  a reproducing kernel Hilbert space.

\subsection{The case of MOT}
For $d=1$, under the (martingale)  twist condition $\partial_{s_1} \partial_{s_2}^2 c \geq 0$, the optimal probability measure $\PP^*$ is shown to be supported not on a single map $T$ but on two maps $T_d(x) \leq x \leq T_u(x)$ \cite{bei2,phl3}:
\beaa \PP^*(s_1,s_2)=q(s_1) \delta(s_2-T_u(s_1)) +(1-q(s_1))\delta(s_2-T_d(s_1)) \eeaa   This leads to the following  minimax optimization:
\beaa {\mathrm{MK}}_c(\mu^1,\mu^2):=\inf_{u \in \mathrm{L}^1(\mu^2), h \in C^0(\RR^d,[0,1]) } \sup_{T_u:\RR \mapsto \RR, T_d:\RR \mapsto \RR,q:\RR \mapsto [0,1]}  \\ \EE^{\mu^1}[ q(S_1) (c(S_1,T_u(S_1))-u(T_u(S_1))-h(S_1)(T_u(S_1)-S_1)) \\ + (1-q(S_1)) (c(S_1,T_d(S_1))-u(T_d(S_1))-h(S_1)(T_d(S_1)-S_1)) ]+\EE^{\mu^2}[u(S_2)] \eeaa
\no Note that the martingale condition leads explicitly to  $q(x):={x-T_d(x) \over T_u(x)-T_d(x)}$ but we do not use this equation in order to preserve the concavity-convexity property with respect to the neural network weights (in the case of one layer). The algorithm is then similar to the one presented for OT except that now we have five (instead of two) neural networks for the potentials $h,u,q$ and the two maps $T_u$ and $T_d$.

\no For $d\geq 2$, one can characterize the cost functions for which the optimal probability measure $\PP^*$ is supported on $n$ maps $T_i$ \cite{demarch2}. The above optimization becomes therefore:
\beaa {\mathrm{MK}}_c(\mu^1,\mu^2):=\inf_{u \in \mathrm{L}^1(\mu^2), h \in C^0(\RR^d,[0,1]) } \sup_{(T_i)_{1 \leq i \leq n}:\RR^d \mapsto \RR^d, q_i:\RR^d \mapsto [0,1]} \\ \EE^{\mu^1}[ q_i(S_1) (c(S_1,T_i(S_1))-u(T_i(S_1))-h(S_1)(T_i(S_1)-S_1))]
+\EE^{\mu^2}[u(S_2)] \eeaa where $q_n:=1-\sum_{i=1}^{n-1} q_i$. In practise, the number of maps $n$ can be seen as an hyperparameter that can be optimized.

\section{Numerical examples} \label{Numerical examples}

\subsection{OT in $d=1$}
\begin{figure}[h]
\begin{center}
\includegraphics[width=7cm,height=5cm]{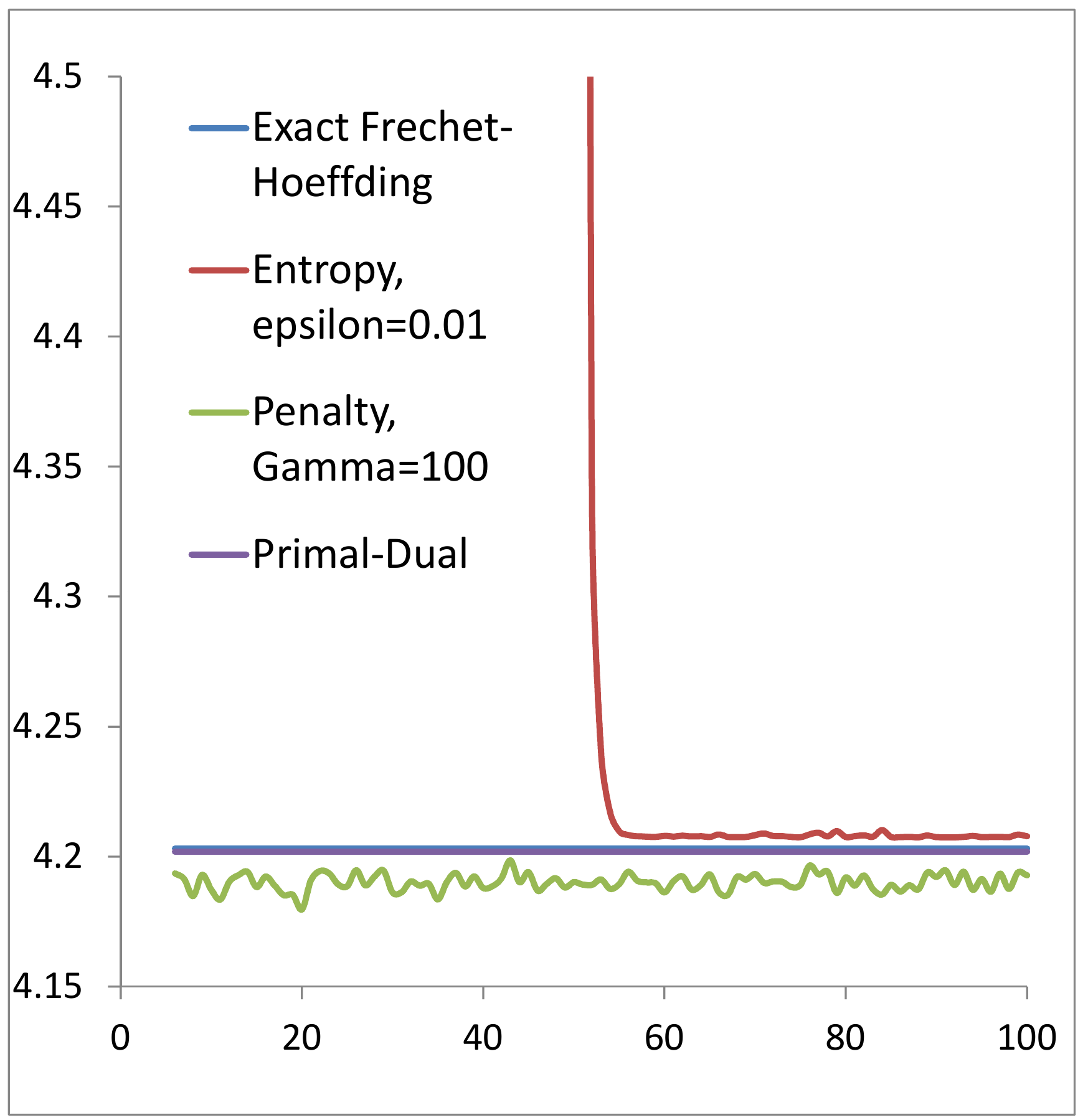}
\includegraphics[width=7cm,height=5cm]{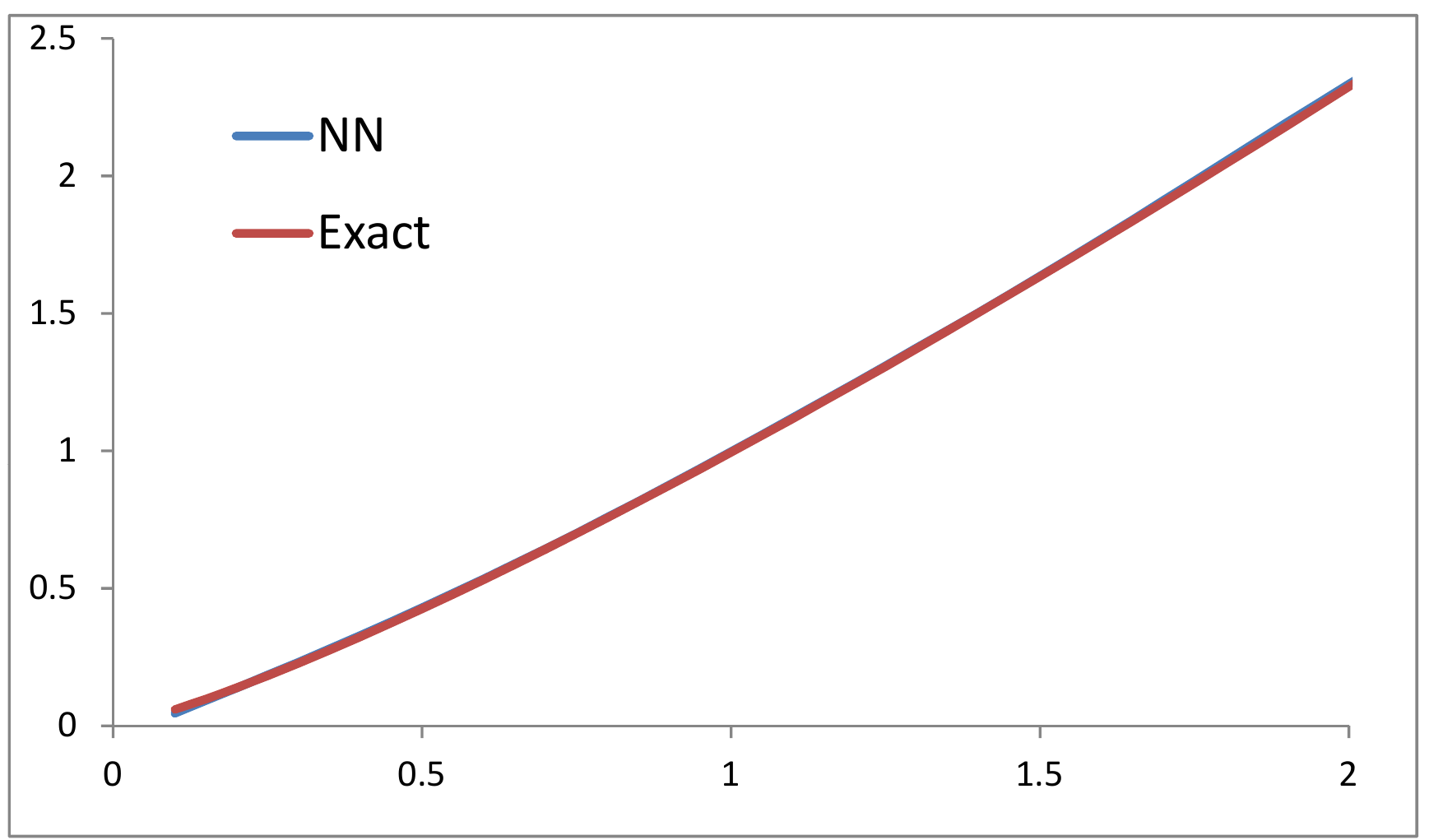}

\end{center}
\caption{OT: $c(s_1,s_2)=(s_1+s_2)^2$. $\mu_1$ and $\mu_2$ are two  log-normal distributions with variances $0.2^2$ and $0.2^2 \times 1.5$. Exact = $4.20$. For the penalization method, we have chosen $\gamma=100$ (similar results for $\gamma=50,200$). The number of iterations has been divided by $10^4$. }
 \label{exa1}
\end{figure}

\begin{figure}[h]
\begin{center}
\includegraphics[width=7cm,height=5cm]{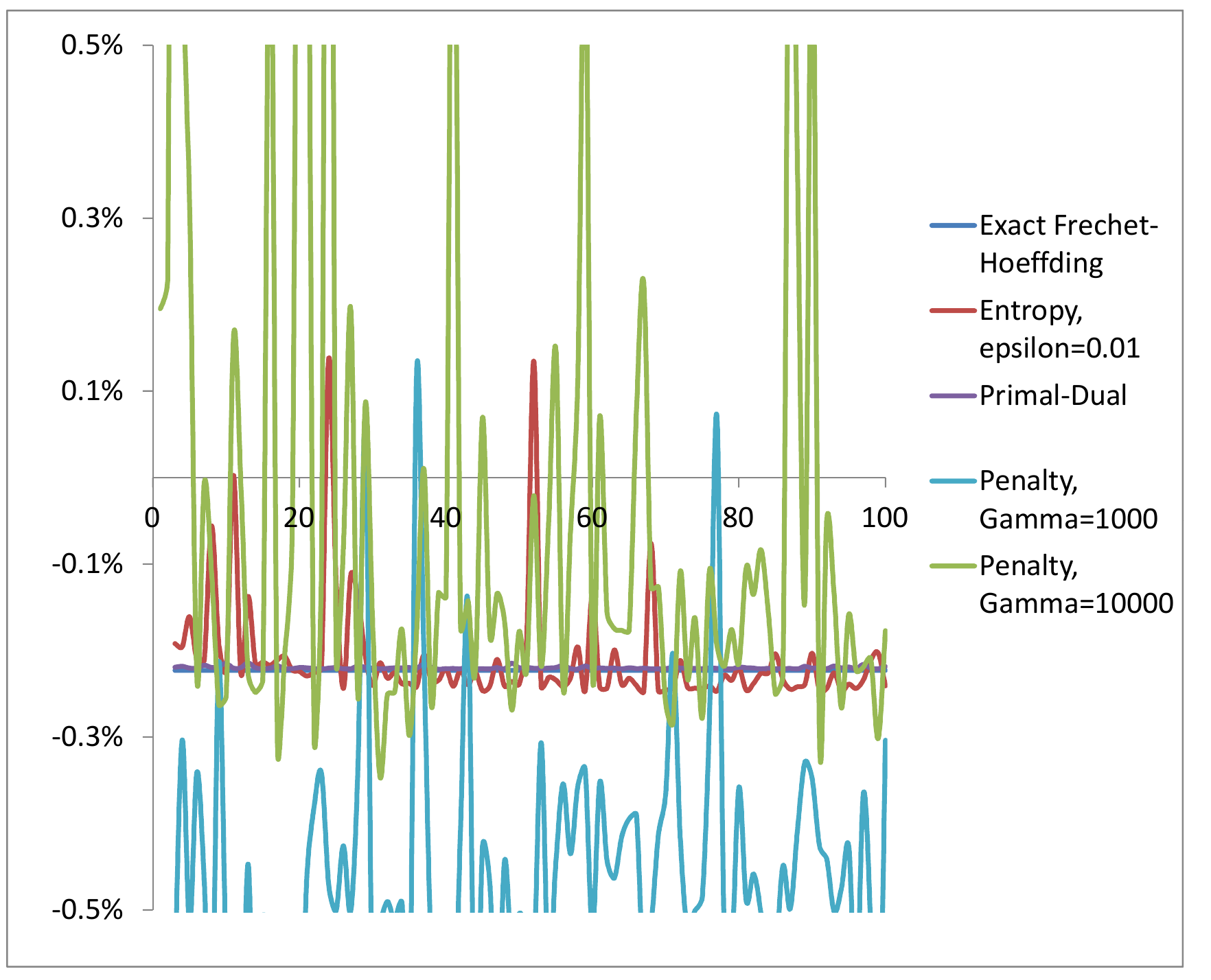}
\includegraphics[width=7cm,height=5cm]{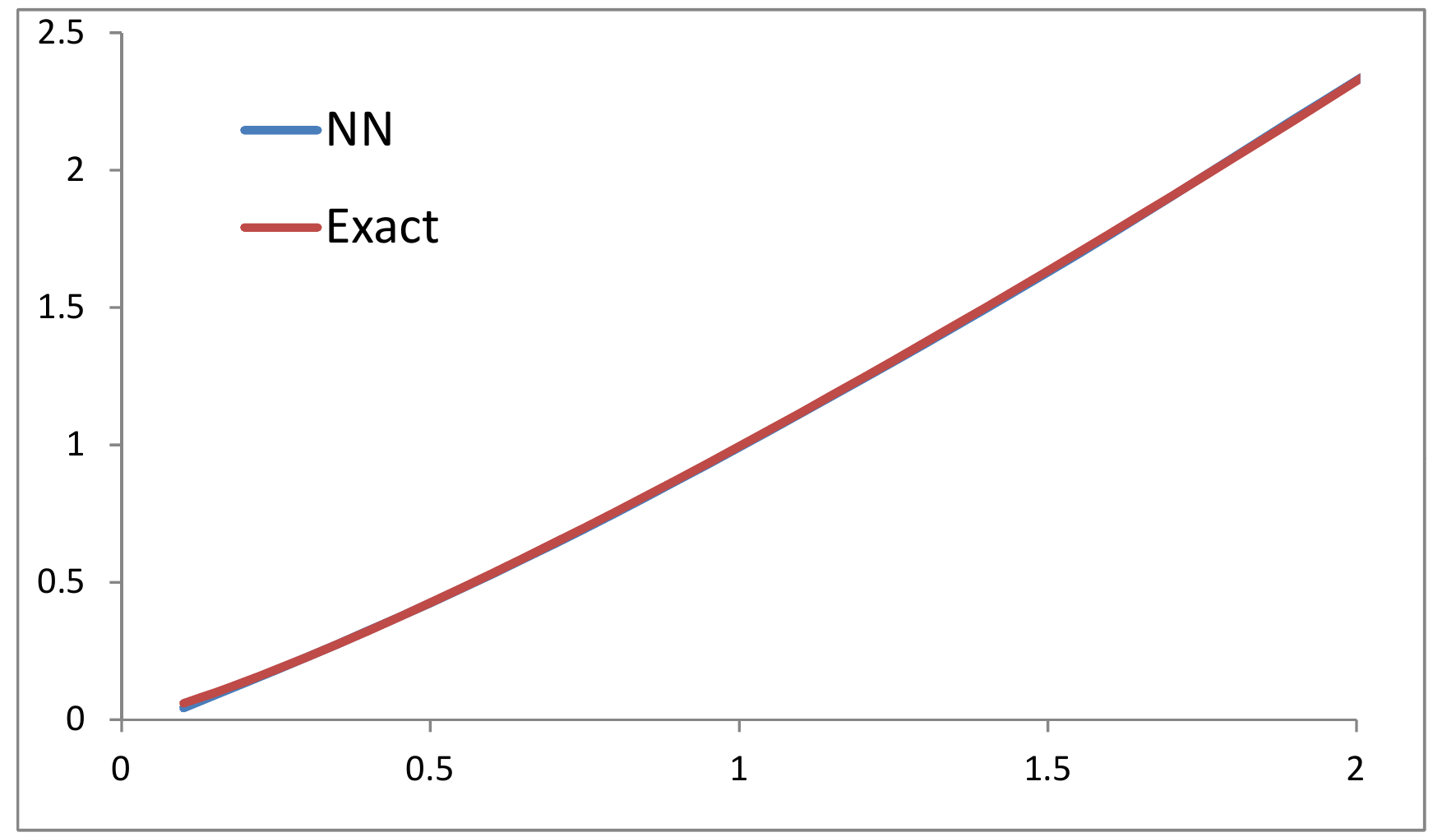}
\end{center}
\caption{OT: $c(s_1,s_2)=-(s_1-s_2)^2$. $\mu_1$ and $\mu_2$ are two  log-normal distributions with variances $0.2^2$ and $0.2^2 \times 1.5$.  Exact = $-0.22\%$. For the penalization method, we have chosen $\gamma=1000,10000$. The number of iterations has been divided by $10^4$.  }
\label{exa2}
\end{figure}

\no We first check our algorithm described in Section \ref{Arrow-Hurwicz algorithm: recipe}  for OT problem in $d=1$. We consider the two cost functions  $c(s_1,s_2)=(s_1+s_2)^2$ and $c(s_1,s_2)=-(s_1-s_2)^2$ satisfying the conditions in Proposition \ref{Convergence} (see Figures \ref{exa1} and \ref{exa2}). $\mu_1$ and $\mu_2$ are chosen to be  two log-normal distributions in $\RR^+$ centered at $S_0=1$ and with variances $0.2^2$ and $0.2^2 \times 1.5$. They are simulated using $2^{13}$ Monte-Carlo paths. For each neural network, we have used $2$ hidden layers of dimension $4$. We have also used a Adam stochastic gradient descent \cite{adam} with $64$ minibatches for the computation of the online gradients and our algorithm has been written from crash in \texttt{C++}. The exact solution has been computed using formula (\ref{HF}) and performing a 1d numerical integration. We have compared our algorithm with the entropy relaxation   and  the penalization methods outlined in Sections \ref{Entropy}-\ref{Penalization}. We can observe that our primal-dual algorithm converges faster (to the exact solution). On one hand, the choice of the gamma factor in the penalization method is tricky. Taking a small value of $\gamma$ results into convergence towards a false solution and a large $\gamma$ gives noisy results. On the other hand,  the entropy relaxation  needs more iterations to converge. We  have used in all our numerical experiments at most $10^6$ iterations. For each $10^4 \times n$ iterations where $n$ ranges from $1$ up to $10^2$, we have computed the functional $J(\theta_n, \omega_n)$ by averaging over our recorded $2^{13}$ Monte-Carlo paths. We have also plotted the map found by our algorithm (denoted ``NN'') and compared with the Fr\'echet-Hoeffding solution $T(s)=F_2^{-1} \circ F_1(s)$. We found a perfect match (the blue and red curves coincide).

\subsection{$2$-Wassertein distance in $\RR^d$, $d=2,10,20$}

\begin{figure}[h]
\begin{center}
\includegraphics[width=7cm,height=5cm]{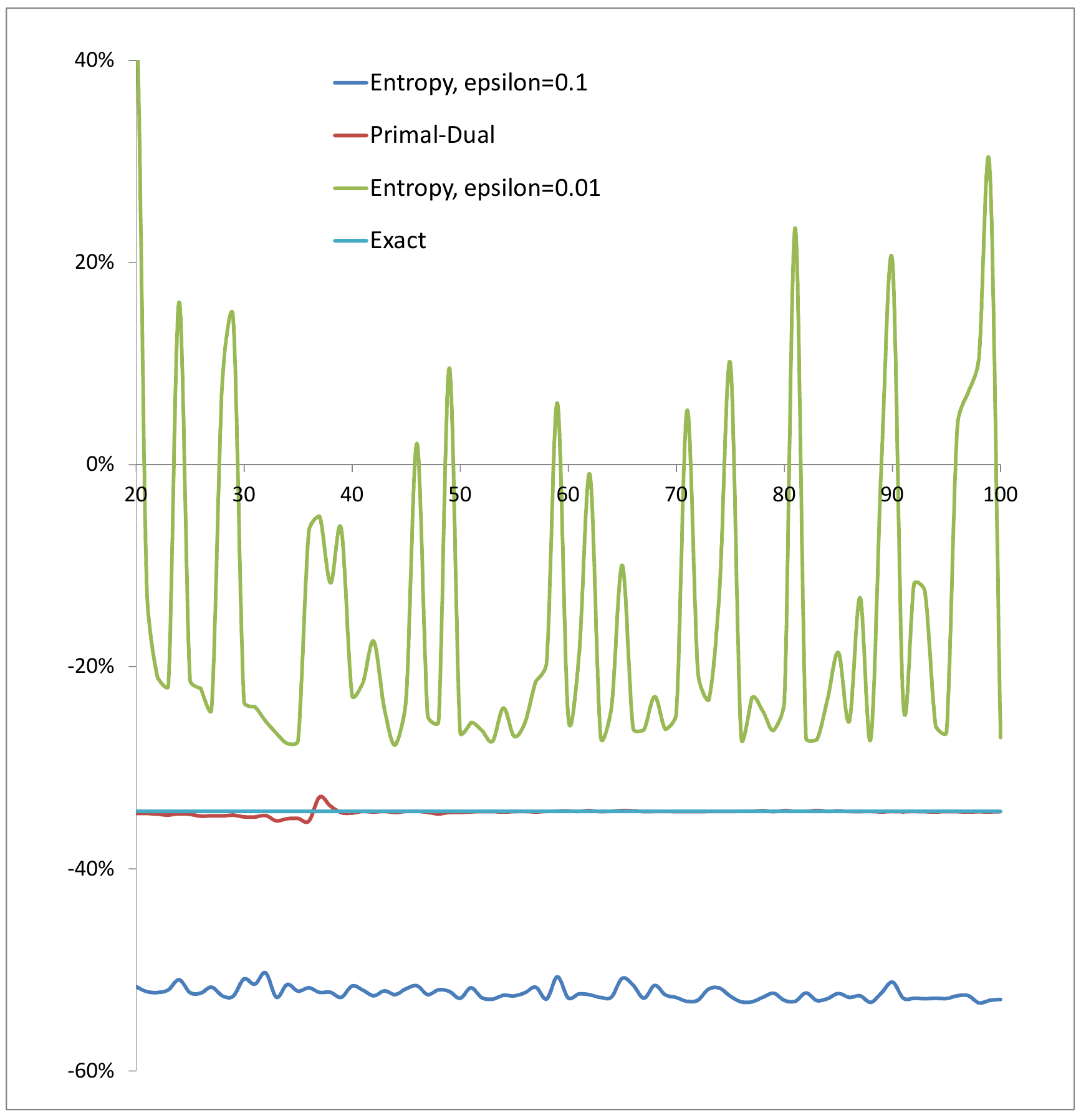}
\includegraphics[width=7cm,height=5cm]{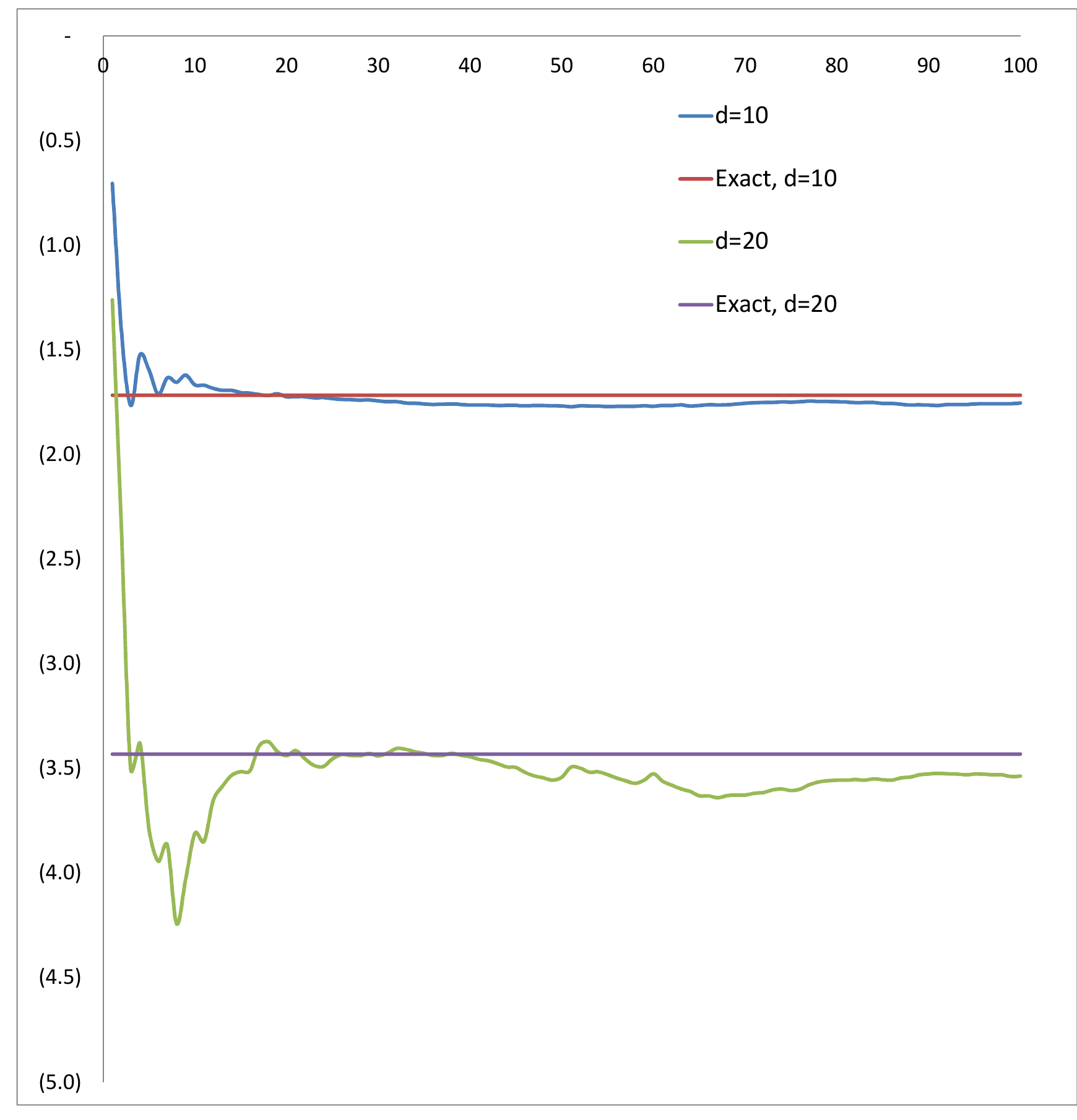}
\end{center}
\caption{OT: $c(s_1,s_2)=-\sum_{i=1}^d (s^i_1-s^i_2)^2$. $\mu_1$ and $\mu_2$ are two uncorrelated normal distributions in $\RR^d$ with variances $1$ and $2$.
\no Left: $d=2$. The number of iterations has been divided by $10^4$. Right: $d=10$ and $d=20$. The number of iterations has been divided by $10^3$ here as our algorithm converges quickly.}
\label{exa3}
\end{figure}

\no Next, we compute the $2$-Wassertein distance in $\RR^d$. In our notation, this corresponds to the payoff $c(s_1,s_2)=-\sum_{i=1}^d (s^i_1-s^i_2)^2$ with a minus sign. We have first considered $d=2$ (see Figure \ref{exa3}--left).  We have compared the entropy relaxation method against our primal-dual algorithm. As concluded in $d=1$, our algorithm converges faster and the entropy relaxation method is unstable according to our choice of $\epsilon$. For large epsilon, the Wasserstein distance is underestimated and for  small epsilon, our SGD is noisy and therefore the result can not be trusted. As a consequence, the entropy relaxation method could not be used as presented for computing the Wasserstein distance.  The convergence is very fast for our primal-dual method. Here $\mu_1$ and $\mu_2$ are chosen to be  two uncorrelated {\it normal} distributions in $\RR^d$ with variances $1$ and $2$ for which the   exact $2$-Wassertein distance in $\RR^d$ is
${\cal W}_2(\mu^1,\mu^2)^2=d (\sqrt{2}-\sqrt{1})^2$. Then, we consider only our primal-dual algorithm and  take  $d=10$ and $d=20$ (see Figure \ref{exa3}--right). For each neural network, we have used $1$ hidden layer of dimension $50$.

\subsection{MOT in $d=1$}
\begin{figure}[h]
\begin{center}
\includegraphics[width=7cm,height=5cm]{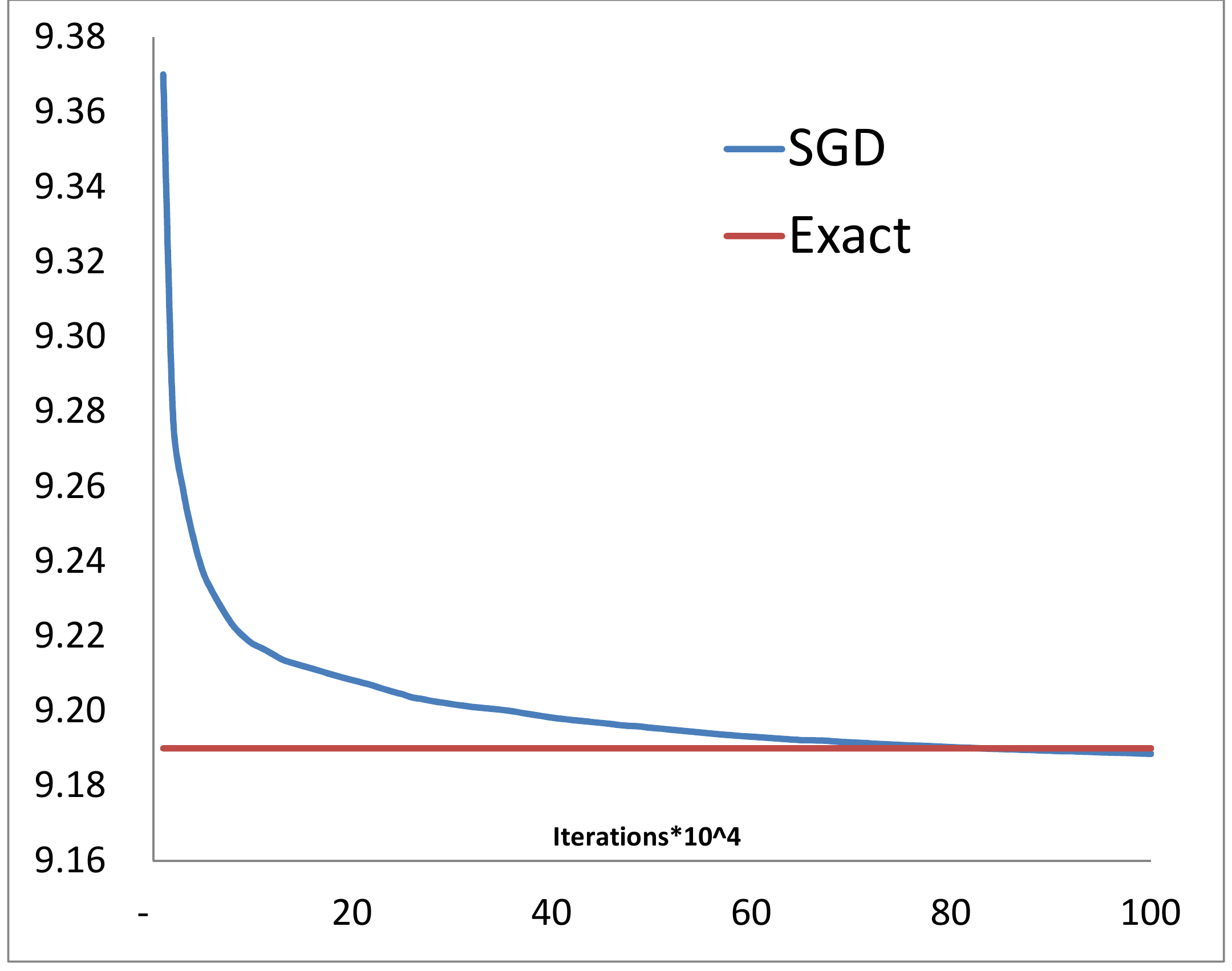}
\end{center}
\caption{MOT: $c(s_1,s_2)=(s_1+s_2)^3$. $\mu_1$ and $\mu_2$ are two  log-normal distributions with variances $0.2^2$ and $0.2^2 \times 1.5$ (in the convex order). Exact using a simplex: $9.19$.}
 \label{exa4}
\end{figure}

\no A similar test has been performed in the case of MOT in $d=1$ with a cost $c(s_1,s_2)=(s_1+s_2)^3$ for which the martingale twist condition $\partial_{s_1}\partial_{s_2}^2 c>0$ is satisfied. Our optimization converges towards the exact solution obtained using a simplex algorithm (see Figure \ref{exa4}).

\subsection{Anomaly detection in $d=2$} 
\begin{figure}[h]
\begin{center}
\includegraphics[width=7cm,height=5cm]{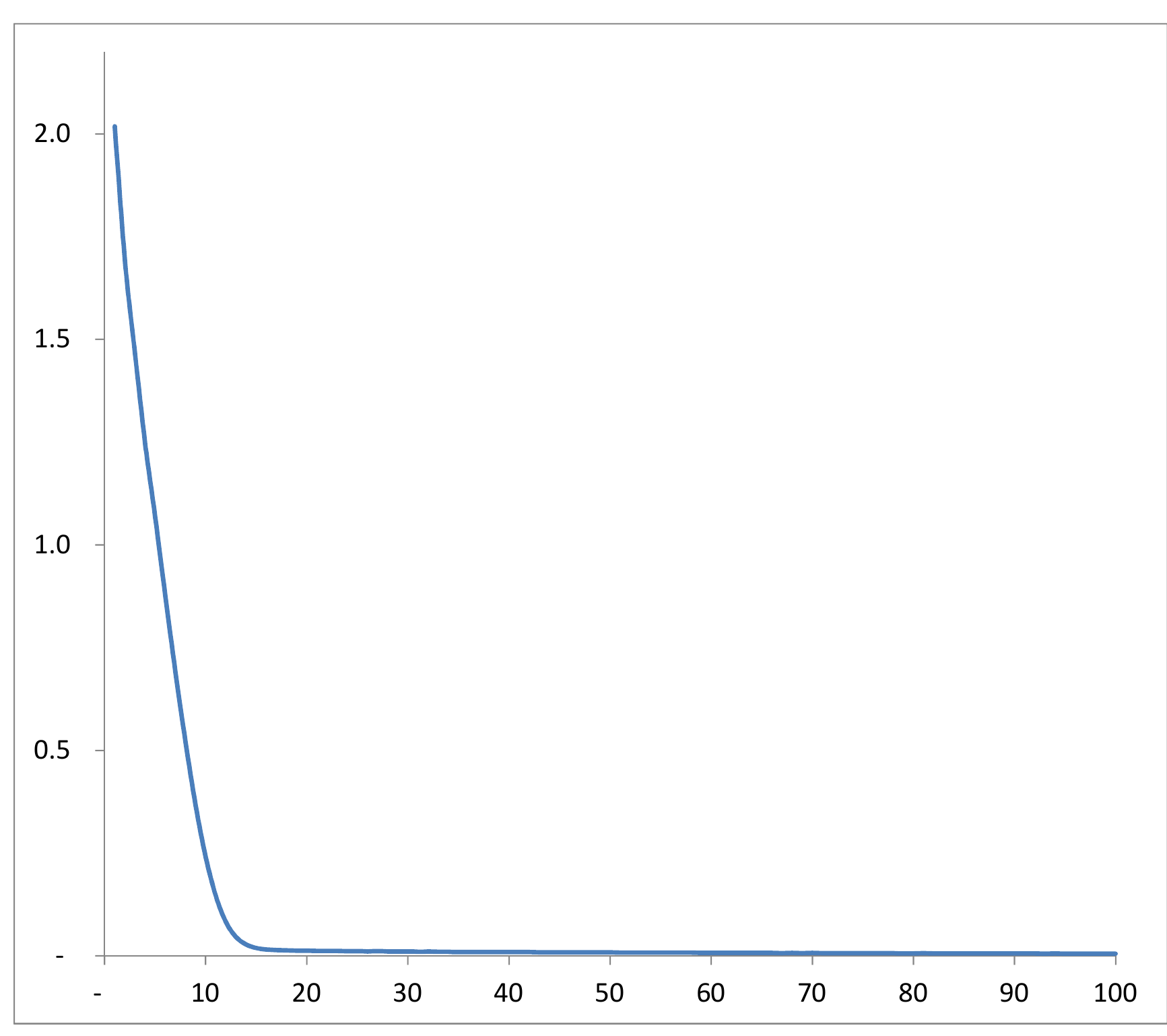}
\end{center}
\caption{Convergence of the $2$-Wasserstein distance ${\cal W}_2(\mu^\mathrm{real},\hat{T}_{\hat{\omega}} \# \mu^0)$. The number of iterations has been divided by $10^4$.}
 \label{exa5}
\end{figure}

\begin{figure}[h]
\begin{center}
\includegraphics[width=7cm,height=5cm]{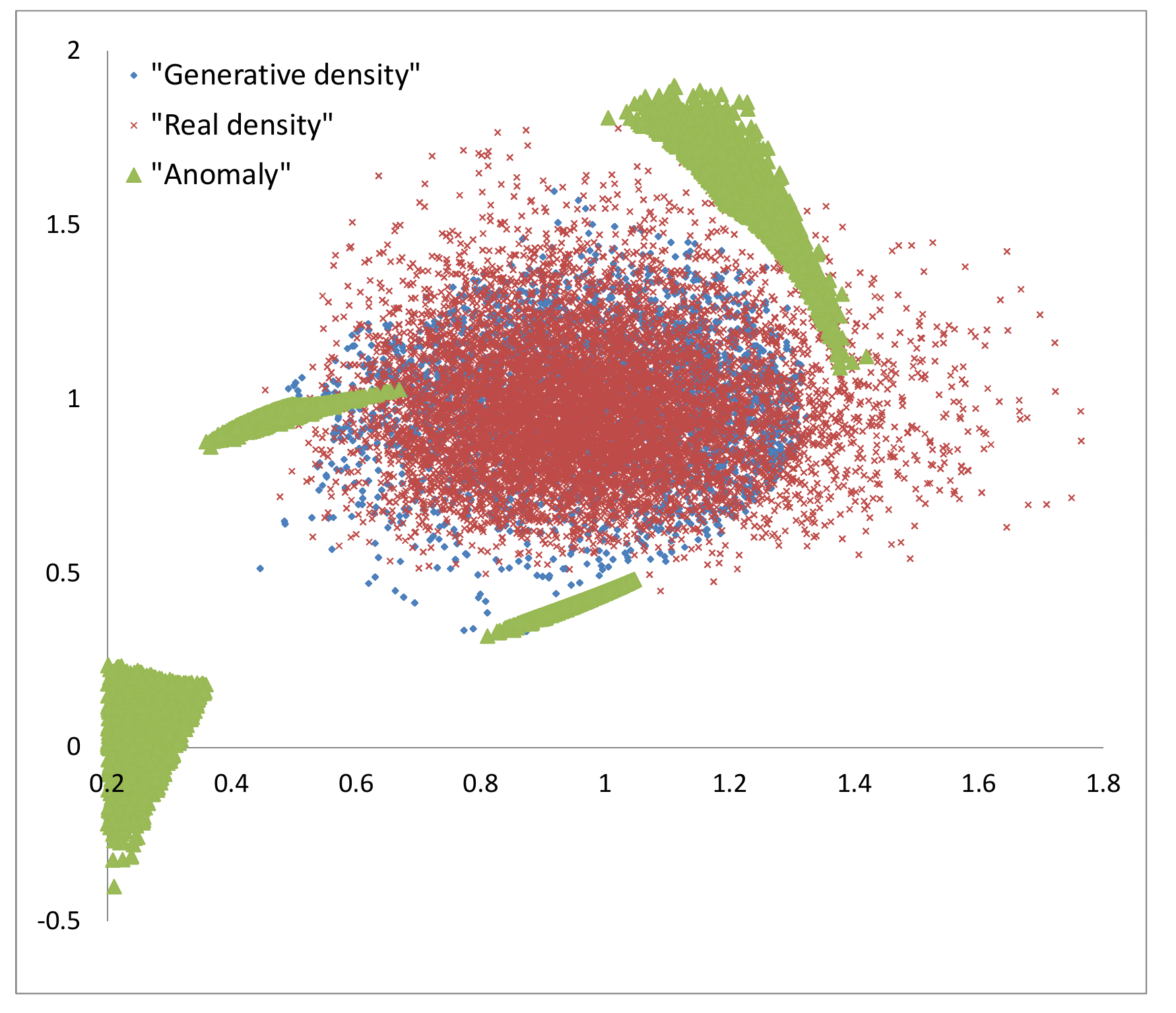}
\end{center}
\caption{Scatter plot of $\hat{T}_{\hat{\omega}}(G)$ with $G \in \mathrm{N}(0,I_2)$ in blue and $\mu^\mathrm{real}$ in red.  Scatter plot of $\hat{T}_{\hat{\omega}}(G+3 \times \mathrm{sign}(G))$ with $G \in \mathrm{N}(0,I_2)$  in green. $\mu^\mathrm{real}$ is a two-dimensional  uncorrelated log-normal distribution with mean $-0.02$, variance $0.04$ and $\mu^\mathrm{0}$ is  a two-dimensional  uncorrelated normal distribution.}
 \label{exa6}
\end{figure}

As a  final simple numerical example, we consider our anomaly
detection algorithm outlined in Section \ref{Anomaly detection and data generator}. We have used $2$ hidden layers of dimension $10$ with linear activation output. We take for $\mu^\mathrm{real}$ a two-dimensional  uncorrelated log-normal distribution with mean $-0.02$, variance $0.04$ and for $\mu^\mathrm{0}$  a two-dimensional  uncorrelated normal distribution. They are simulated using $2^{13}$ Monte-Carlo paths. Note that the stochastic Arrow-Hurwicz iterations over $u_\theta$ and $T_\omega$ are performed and each $1000$ iterations, a stochastic gradient descent minimization over $\hat{T}_{\hat{\omega}}$ is done. We have plotted in Figure \ref{exa5} the $2$-Wasserstein distance ${\cal W}_2(\mu^\mathrm{real},{\hat{T}_{\hat{\omega}}} \# \mu^0)$ each $10^4$ iterations and this converges, as expected, to zero. Once the mapping $\hat{T}: \RR^2 \mapsto \RR^2$ is constructed by optimization, we generate some ``anomalies''  $\hat{T}_{\hat{\omega}}(G+3 \times \mathrm{sign}(G))$ by drawing some normal  variables $ G \in \mathrm{N}(0,I_2)$ in $\RR^2$ and adding  an anomaly factor $3\times \mathrm{sign}(G)$.  The ``normal'' variables  $\hat{T}_{\hat{\omega}}(G)$ are generated without introducing this anomaly factor. The two-dimensional ``normal'' and ``abnormal'' variables generated are then displayed in Figure \ref{exa6}. As expected, the ``abnormal'' data live on the edge of the  two-dimensional  uncorrelated log-normal distribution $\mu^\mathrm{real}$, which is close to $\hat{T}_\#\mu_0$ with respect to the $2$-Wasserstein distance (see  Figure \ref{exa5}).

\end{document}